\documentclass[amstex,12pt,russian,amssymb]{article}

\usepackage{mathtext}
\usepackage[cp1251]{inputenc}
\usepackage[T2A]{fontenc}
\usepackage[russian]{babel}
\usepackage[dvips]{graphicx}
\usepackage{amsmath}
\usepackage{amssymb}
\usepackage{amsxtra}
\usepackage{latexsym}
\usepackage{ifthen}

\begin{document}

\newcounter{lemma}
\newcommand{\lemma}{\par \refstepcounter{lemma}%
{\bf Лемма \arabic{lemma}.}}

\newcounter{corollary}
\newcommand{\corollary}{\par \refstepcounter{corollary}%
{\bf Следствие \arabic{corollary}.}}

\newcounter{remark}
\newcommand{\remark}{\par \refstepcounter{remark}%
{\bf Замечание \arabic{remark}.}}

\newcounter{theorem}
\newcommand{\theorem}{\par \refstepcounter{theorem}%
{\bf Теорема \arabic{theorem}.}}

\newcounter{proposition}
\newcommand{\proposition}{\par \refstepcounter{proposition}%
{\bf Предложение \arabic{proposition}.}}

\renewcommand{\refname}{\centerline{\bf Список литературы}}

\newcommand{\proof}{{\it Доказательство.\,\,}}

\noindent УДК 517.5

{\bf Е.А. Севостьянов} (Институт прикладной математики и механики
НАН Украины, Донецк, Украина)

{\bf Є.О. Севостьянов} (Інститут прикладної математики і механіки
НАН України, Донецьк, Україна)

{\bf E.A. Sevost'yanov} (Institute of Applied Mathematics and
Mechanics of NAS of Ukraine, Donetsk, Ukraine)

{\bf О теоремах сходимости и нормальности локальных гомеоморфизмов
классов Орлича--Соболева }

{\bf Про теореми збіжності та нормальності локальних гомеоморфізмів
класів Орліча--Соболєва}

{\bf On convergence and normality theorems of local homeomorphisms
of Orlicz--Sobolev classes}

\medskip
Исследуются проблемы, связанные с локальными и предельными
свойствами классов Орлича--Соболева конечного искажения, активно
изучаемых в последнее время. Показано, что локально равномерным
пределом локальных гомеоморфизмов классов Орлича--Соболева является
локальный гомеоморфизм, либо постоянная, как только их дилатации
удовлетворяют соответствующим ограничениям. Кроме того, при
некоторых дополнительных условиях показано, что локально равномерным
пределом отображений указанного класса может быть только локальный
гомеоморфизм либо постоянная.

\medskip
Досліджено проблеми, пов'язані з локальними і межовими властивостями
класів Орліча--Соболєва скінченного спотворення, які активно
вивчаються останнім часом. Показано, що локально рівномірною
границею локальних гомеоморфізмів класів Орліча--Соболєва є
локальний гомеоморфізм, або стала, як тільки їх дилатації
задовольняють відповідні обмеження. Крім того, за деякими
додатковими умовами показано, що локально рівномірною границею
відображень вказаного класу може бути тільки локальний гомеоморфізм
або стала.

\medskip
There are investigated problems connected with local and boundary
properties of Orlicz--Sobolev classes of finite distortion which are
actively studied last time. It is showed that, a locally uniform
limit of local homeomorphisms of Orlicz--Sobolev class is a local
ho\-me\-o\-mor\-phism or a constant whenever it's dilatations
satisfy corresponding restrictions. Besides that, at some additional
conditions, it is showed that, a locally uniform limit of the
mappings mentioned above is a local ho\-me\-o\-mor\-phism or a
constant.

\newpage

{\bf 1. Введение.} Настоящая заметка посвящена изучению отображений
с конечным искажением, активно изучаемых в последнее время (см.,
напр., \cite{IM}, \cite{MRSY}--\cite{MRSY$_1$} и \cite{GRSY}). В
частности, речь идёт об изучении классов Орлича--Соболева,
определяющихся как семейства отображений, градиент которых локально
принадлежит некоторому классу Орлича. Точнее, пусть
$\varphi:[0,\infty)\rightarrow[0,\infty)$ -- неубывающая функция,
$f$ -- локально интегрируемая вектор-функция $n$ вещественных
переменных $x_1,\ldots,x_n,$ $f=(f_1,\ldots,f_m),$ $f_i\in
W_{loc}^{1,1},$ $i=1,\ldots,m.$ Будем говорить, что $f:D\rightarrow
{\Bbb R}^n$ принадлежит классу $W^{1,\varphi}_{loc},$ пишем $f\in
W^{1,\varphi}_{loc},$ если
%
$$\int\limits_{G}\varphi\left(|\nabla
f(x)|\right)\,dm(x)<\infty$$
%
для любой компактной подобласти $G\subset D,$ где $|\nabla
f(x)|=\sqrt{\sum\limits_{i=1}^m\sum\limits_{j=1}^n\left(\frac{\partial
f_i}{\partial x_j}\right)^2}.$ Класс $W^{1,\varphi}_{loc}$
называется классом {\it Орлича--Соболева}. Как было показано в одной
из совместных работ автора, семейства гомеоморфизмов класса
Орлича--Соболева с конечным искажением являются нормальными
(равностепенно непрерывными) при определённых дополнительных
условиях на характеристику квазиконформности отображений и
количество выпускаемых этими отображениями значений (см., напр.,
\cite[теоремы~7 и 9, следствие~13]{KRSS}). Кроме того, ранее было
установлено, что локально равномерным пределом указанных отображений
является гомеоморфизм либо постоянная (см. \cite[следствие~9]{KRSS}
и \cite[раздел~4]{RSS}). В настоящей заметке будет показано, что
соответствующие результаты верны также и для аналогичных семейств
локальных гомеоморфизмов. Следует отметить, что исследование
локальных гомеоморфизмов, не являющихся гомеоморфизмами, требует
привлечения существенно иной техники. (Отображение $f:D\rightarrow
{\Bbb R}$ называется здесь и далее локальным гомеоморфизмом, если
каждая точка $x_0\in D$ имеет окрестность $U,$ такую что $f|_{U}$
является гомеоморфизмом, т.е., взаимнооднозначным отображением,
обратное к которому также, как и $f,$ является непрерывным).

\medskip
Сформулируем основные утверждения, которые будут доказаны в
настоящей работе.

\medskip
Всюду далее $D$ -- область в ${\Bbb R}^n,$ $n\ge 2,$ $m$ -- мера
Лебега в ${\Bbb R}^n$ и ${\rm dist\,}(A,B)$ -- евклидово расстояние
между множествами $A$ и $B$ в ${\Bbb R}^n.$ Запись $f:D\rightarrow
{\Bbb R}^n$ предполагает, что отображение $f$ непрерывно в $D.$ В
дальнейшем ${\mathcal H}^k$ -- нормированная $k$-мерная мера
Хаусдорфа в ${\Bbb R}^n,$ $1\le k\le n,$
$$B(x_0, r)=\left\{x\in{\Bbb R}^n: |x-x_0|< r\right\}\,,\quad {\Bbb B}^n
:= B(0, 1)\,,$$ $$S(x_0,r) = \{ x\,\in\,{\Bbb R}^n :
|x-x_0|=r\}\,,\quad{\Bbb S}^{n-1}:=S(0, 1)\,,$$
$$A(r_1,r_2,x_0)=\{ x\,\in\,{\Bbb R}^n : r_1<|x-x_0|<r_2\}\,,$$
$\omega_{n-1}$ обозначает площадь единичной сферы ${\Bbb S}^{n-1}$ в
${\Bbb R}^n,$ $\Omega_{n}$ -- объём единичного шара ${\Bbb B}^{n}$ в
${\Bbb R}^n.$ В дальнейшем $J(x, f)={\rm det}\, f^{\,\prime}(x)$ --
{\it якобиан отображения} $f$ в точке $x,$ где $f^{\,\prime}(x)$ --
{\it матрица Якоби} отображения $f$ в точке $x.$

\medskip
Отображение $f:D\rightarrow {\Bbb R}^n$ называется {\it отображением
с конечным искажением}, пишем $f\in FD,$ если $f\in
W_{loc}^{1,1}(D)$ и для некоторой функции $K(x): D\rightarrow
[1,\infty)$ выполнено условие
$\Vert f^{\,\prime}\left(x\right) \Vert^{n}\le K(x)\cdot |J(x,f)|$
при почти всех $x\in D$ (см. \cite[п.~6.3, гл.~VI]{IM}. Для
отображений с конечным искажением корректно определена и почти всюду
конечна так называемая {\it внешняя дилатация} $K_O(x,f)$
отображения $f$ в точке $x,$ определяемая соотношением
\begin{equation}\label{eq0.1.1A}
K_O(x,f)\quad =\quad\left\{
\begin{array}{rr}
\frac{\Vert f^{\,\prime}(x)\Vert^n}{|J(x, f)|}, & J(x,f)\ne 0,\\
1,  &  f^{\,\prime}(x)=0, \\
\infty, & \text{в\,\,остальных\,\,случаях}
\end{array}
\right.\,.
\end{equation}
Для отображения $f:D\,\rightarrow\,{\Bbb R}^n,$ множества $E\subset
D$ и $y\,\in\,{\Bbb R}^n,$  определим {\it функцию кратности $N(y,
f, E)$} как число прообразов точки $y$ во множестве $E,$ т.е.
\begin{equation}\label{eq1.7A}
N(y, f, E)\,=\,{\rm card}\,\left\{x\in E: f(x)=y\right\}\,,
%
N(f, E)\,=\,\sup\limits_{y\in{\Bbb R}^n}\,N(y, f, E)\,.
\end{equation}

\medskip
Пусть $\left(X,\,d\right)$ и
$\left(X^{\,{\prime}},{d}^{\,{\prime}}\right)$ -- метрические
пространства с расстояниями  $d$  и ${d}^{\,{\prime}}$
соответственно. Семейство $\frak{F}$ непрерывных отображений
$f:X\rightarrow {X}^{\,\,\prime}$ называется {\it нормальным}, если
из любой последовательности отображений $f_{m} \in \frak{F}$ можно
выделить подпоследовательность $f_{m_{k}}$, которая сходится
локально равномерно в $X$ к непрерывной функции
$f:\,X\,\rightarrow\, X^{\,\prime}.$

\medskip
Введенное  понятие  очень  тесно  связано  со  следующим. Семейство
$\frak{F}$ отображений $f:X\rightarrow {X}^{\,\prime}$ называется
{\it равностепенно непрерывным в точке} $x_0 \in X,$ если для любого
$\varepsilon>0$ найдётся $\delta>0$ такое, что ${d}^{\,\prime}
\left(f(x),f(x_0)\right)<\varepsilon$ для всех $x$ таких, что
$d(x,x_0)<\delta$ и для всех $f\in \frak{F}.$ Говорят, что
$\frak{F}$ {\it равностепенно непрерывно}, если $\frak{F}$
равностепенно непрерывно в каждой  точке $x_0\in X.$ Согласно одной
из версий теоремы Арцела-Асколи (см., напр., \cite[пункт~20.4]{Va}),
если $\left(X,\,d\right)$ -- сепарабельное метрическое пространство,
а $\left(X^{\,\prime},\, d^{\,\prime}\right)$ -- компактное
метрическое пространство, то семейство $\frak{F}$ отображений
$f:X\rightarrow {X}^{\,\prime}$ нормально тогда  и только тогда,
когда  $\frak{F}$ равностепенно непрерывно.

\medskip
Отметим, что всюду далее, если не оговорено противное, $(X, d)=(D,
|\cdot|),$ где $D$ -- область в ${\Bbb R}^n,$ а $|\cdot|$ --
евклидова метрика, $|x-y|=\sqrt{\sum\limits_{i=1}^n (y_i-x_i)^2},$
$x=(x_1,\ldots,x_n),$ $y=(y_1,\ldots,y_n);$ $\left(X^{\,\prime},\,
d^{\,\prime}\right)=\left(\overline{{\Bbb R}^n},\, h\right),$ где
$\overline{{\Bbb R}^n}={\Bbb R}^n\cup\{\infty\},$ $h$ -- хордальная
метрика,
$$h(x,\infty)=\frac{1}{\sqrt{1+{|x|}^2}}\,,\quad h(x,y)=\frac{|x-y|}{\sqrt{1+{|x|}^2} \sqrt{1+{|y|}^2}}\,,  x\ne
\infty\ne y\,.$$
При заданном множестве $E\subset\overline{{\Bbb R}^n}$ полагаем
$h(E)=\sup\limits_{x, y\in E}h(x, y)$ -- хордальный диаметр
множества $E.$

\medskip
Для заданного числа $\delta>0,$ неубывающей функции
$\varphi:[0,\infty)\rightarrow[0,\infty),$ измеримой по Лебегу
функции $Q:D\rightarrow [1, \infty]$ и числа $N\in {\Bbb N}$
обозначим символом $\frak{R}_{\varphi, Q, N, \delta}$ семейство всех
локальных гомеоморфизмов $f:D\rightarrow {\Bbb R}^n$ класса
$W^{1,\varphi}_{loc},$ имеющих конечное искажение, таких что $N(f ,
D)\le N,$ $K_O^{n-1}(x, f)\le Q(x)$ и $h(\overline{{\Bbb
R}^n}\setminus f(D))\ge \delta.$ Справедлива следующая

\medskip
\begin{theorem}\label{th1}
{\sl\, Пусть $n\ge 3,$ тогда семейство отображений
$\frak{R}_{\varphi, Q, N, \delta}$ является равностепенно
непрерывным в некоторой фиксированной точке $x_0\in D,$ если $Q\in
L_{loc}^1(D),$
\begin{equation}\label{eqOS3.0a}
\int\limits_{1}^{\infty}\left[\frac{t}{\varphi(t)}\right]^
{\frac{1}{n-2}}dt<\infty
\end{equation}
и, кроме того, при некотором $\varepsilon_0>0,$ $\varepsilon_0<{\rm
dist}(x_0, \partial D),$ выполнено следующее условие расходимости
интеграла:
\begin{equation}\label{eq9}
\int\limits_{0}^{\varepsilon_0}
\frac{dt}{tq_{x_0}^{\,\frac{1}{n-1}}(t)}=\infty,
\end{equation}
где, как обычно,
$q_{x_0}(r):=\frac{1}{\omega_{n-1}r^{n-1}}\int\limits_{|x-x_0|=r}Q(x)\,d{\mathcal
H}^{n-1}$ -- среднее интегральное значение функции $Q$ над сферой
$S(x_0, r).$ В частности, заключение теоремы \ref{th1} является
верным, если $q_{x_0}(r)=\,O\left({\left[
\log{\frac{1}{r}}\right]}^{n-1}\right)$ при $r\rightarrow 0.$}
\end{theorem}

\medskip
Из приведённого выше критерия Арцела-Асколи вытекает следующее

\medskip
\begin{corollary}\label{cor2}{\sl\,
В условиях теоремы \ref{th1} семейство отображений
$\frak{R}_{\varphi, Q, N, \delta}$ является нормальным семейством
отображений, как только условие (\ref{eq9}) выполнено в каждой точке
$x_0$ области $D.$}
\end{corollary}

\medskip
Сформулируем ещё один важнейший результат работы.

\medskip
Будем говорить, что локально интегрируемая функция
${\varphi}:D\rightarrow{\Bbb R}$ имеет {\it конечное среднее
колебание} в точке $x_0\in D$, пишем $\varphi\in FMO(x_0),$ если
%
%
%
%
$\limsup\limits_{\varepsilon\rightarrow
0}\frac{1}{\Omega_n\varepsilon^n}\int\limits_{B( x_0,\,\varepsilon)}
|{\varphi}(x)-\overline{{\varphi}}_{\varepsilon}|\, dm(x)<\infty,$
%
%
где
$\overline{{\varphi}}_{\varepsilon}=\frac{1}
{\Omega_n\varepsilon^n}\int\limits_{B(x_0,\,\varepsilon)}
{\varphi}(x)\, dm(x).$
\medskip
Заметим, что, как известно, $\Omega_n\varepsilon^n=m(B(x_0,
\varepsilon)).$ Имеет место следующая

\medskip
\begin{theorem}\label{th2}
{\sl\, При $n\ge 3$ семейство отображений $\frak{R}_{\varphi, Q, N,
\delta}$ является равностепенно непрерывным в точке $x_0\in D,$ если
выполнено условие (\ref{eqOS3.0a}) и, кроме того, $Q\in FMO(x_0).$}
\end{theorem}

\medskip
Из теоремы \ref{th2} на основании приведённого выше критерия
Арцела-Асколи вытекает следующее

\medskip
\begin{corollary}\label{cor3}{\sl\,
В условиях теоремы \ref{th2} семейство отображений
$\frak{R}_{\varphi, Q, N, \delta}$ является нормальным семейством
отображений, как только условие $Q\in FMO(x_0)$ выполнено в каждой
точке $x_0$ области $D.$}
\end{corollary}

\medskip
Ещё одно важное утверждение, которое будет доказано в настоящей
работе, может быть сформулировано следующим образом.

\medskip
\begin{theorem}\label{th3}
{\sl Пусть $n\ge 3,$ $f_m$ -- последовательность локальных
гомеоморфизмов класса $W^{1,\varphi}_{loc},$ имеющих конечное
искажение, таких что $N(f , D)\le N,$ $K_O^{n-1}(x, f)\le Q(x),$ где
$\varphi:[0,\infty)\rightarrow[0,\infty)$ -- неубывающая функция,
$Q:D\rightarrow [1, \infty]$ -- измеримая по Лебегу функция, а $N\in
{\Bbb N}$ -- некоторое фиксированное число. Если последовательность
$f_m$ сходится локально равномерно по метрике $h$ к отображению
$f:D\rightarrow \overline{{\Bbb R}^n},$ функция $\varphi$
удовлетворяет условию (\ref{eqOS3.0a}), а функция $Q$ в каждой точке
$x_0\in D$ удовлетворяет одному из условий: (\ref{eq9}) либо $Q\in
FMO(x_0),$ то тогда либо $f$ -- локальный гомеоморфизм из $D$ в
${\Bbb R}^n,$ либо $f$ -- постоянная из $D$ в $\overline{{\Bbb
R}^n}.$}
\end{theorem}

\medskip
{\bf 2. Вспомогательные результаты}. Следует отметить, что
доказательство основных результатов работы опирается на некоторый
аппарат, на первый взгляд, не связанный с классами Орлича и
Орлича--Соболева. Речь идёт о так называемых нижних $Q$-отображениях
и кольцевых $Q$-отображениях (см., напр., \cite{MRSY}). Напомним
некоторые определения, связанные с понятием поверхности, интеграла
по поверхности, а также модулей семейств кривых и поверхностей.

\medskip Пусть $\omega$ -- открытое множество в $\overline{{\Bbb
R}^k}:={\Bbb R}^k\cup\{\infty\},$ $k=1,\ldots,n-1.$ Непрерывное
отображение $S:\omega\rightarrow{\Bbb R}^n$ будем называть {\it
$k$-мерной поверхностью} $S$ в ${\Bbb R}^n.$ Число прообразов
$$N(y, S)={\rm card}\,S^{-1}(y)={\rm card}\,\{x\in\omega:S(x)=y\},\
y\in{\Bbb R}^n$$ будем называть {\it функцией кратности} поверхности
$S.$ Другими словами, $N(y, S)$ -- кратность накрытия точки $y$
поверхностью $S.$ Пусть $\rho:{\Bbb R}^n\rightarrow\overline{{\Bbb
R}^+}$ -- борелевская функция, в таком случае интеграл от функции
$\rho$ по поверхности $S$ определяется равенством:  $$\int\limits_S
\rho\,d{\mathcal{A}}:= \int\limits_{{\Bbb R}^n}\rho(y)\,N(y,
S)\,d{\mathcal H}^ky\,.$$
Пусть $\Gamma$ -- семейство $k$-мерных поверхностей $S.$ Борелевскую
функцию $\rho:{\Bbb R}^n\rightarrow\overline{{\Bbb R}^+}$ будем
называть {\it допустимой} для семейства $\Gamma,$ сокр. $\rho\in{\rm
adm}\,\Gamma,$ если
\begin{equation}\label{eq8.2.6}\int\limits_S\rho^k\,d{\mathcal{A}}\ge 1\end{equation}
для каждой поверхности $S\in\Gamma.$
Для заданного числа $p\in(0,\infty)$  {\it $p$-модулем} семейства
$\Gamma$ назовём величину
$$M_p(\Gamma)=\inf_{\rho\in{\rm adm}\,\Gamma} \int\limits_{{\Bbb
R}^n}\rho^p(x)\,dm(x)\,.$$ Мы также полагаем
$M(\Gamma)=M_n(\Gamma),$
а величину $M(\Gamma)$ в этом случае называем {\it модулем
семейства} $\Gamma.$ Заметим, что модуль семейств поверхностей,
определённый таким образом, представляет собой внешнюю меру в
пространстве всех $k$-мерных поверхностей (см. \cite{Fu}).

Пусть $p\ge 1.$ Говорят, что некоторое свойство $P$ выполнено для
{\it $p$-почти всех поверхностей} области $D,$ если оно имеет место
для всех поверхностей, лежащих в $D,$ кроме, быть может, некоторого
их подсемейства, $p$-модуль которого равен нулю. (Как правило, если
речь идёт о конформном модуле, говорят, что указанное свойство
выполнено для {\it почти всех поверхностей} области $D,$ опуская
приставку $"$$n$$"$ в выражении $"$$n$-почти всех$"$). В частности,
говорят, что некоторое свойство выполнено для {\it $p$-почти всех
кривых} области $D$, если оно имеет место для всех кривых, лежащих в
$D$, кроме, быть может, некоторого их подсемейства, $p$-модуль
которого равен нулю.

Будем говорить, что измеримая по Лебегу функция $\rho:{\Bbb
R}^n\rightarrow\overline{{\Bbb R}^+}$ {\it $p$-обобщённо допустима}
для семейства $\Gamma$ $k$-мерных поверхностей $S$ в ${\Bbb R}^n,$
сокр. $\rho\in{\rm ext}_p\,{\rm adm}\,\Gamma,$ если соотношение
(\ref{eq8.2.6}) выполнено для $p$-почти всех поверхностей $S$
семейства $\Gamma.$ {\it Обобщённый $p$-модуль} $\overline
M_p(\Gamma)$ семейства $\Gamma$ определяется равенством
$$\overline M_p(\Gamma)= \inf\int\limits_{{\Bbb
R}^n}\rho^p(x)\,dm(x)\,,$$
где точная нижняя грань берётся по всем функциям $\rho\in{\rm
ext}_p\,{\rm adm}\,\Gamma.$ В случае $p=n$ мы используем обозначения
$\overline M(\Gamma)$ и $\rho\in{\rm ext}\,{\rm adm}\,\Gamma,$
соответственно. Очевидно, что при каждом $p\in(0,\infty),$
$k=1,\ldots,n-1,$ и каждого семейства $k$-мерных поверхностей
$\Gamma$ в ${\Bbb R}^n,$ выполнено равенство $\overline
M_p(\Gamma)=M_p(\Gamma).$

Следующий класс отображений представляет собой обобщение
квазиконформных отображений в смысле кольцевого определения по
Герингу (\cite{Ge$_3$}) и отдельно исследуется различными авторами
(см., напр., \cite[глава~9]{MRSY}). Пусть $D$ и $D^{\,\prime}$ --
заданные области в $\overline{{\Bbb R}^n},$ $n\ge 2,$
$x_0\in\overline{D}\setminus\{\infty\}$ и $Q:D\rightarrow(0,\infty)$
-- измеримая по Лебегу функция. Будем говорить, что $f:D\rightarrow
D^{\,\prime}$ -- {\it нижнее $Q$-отображение в точке} $x_0,$ как
только
\begin{equation}\label{eq1A}
M(f(\Sigma_{\varepsilon}))\ge \inf\limits_{\rho\in{\rm
ext\,adm}\,\Sigma_{\varepsilon}}\int\limits_{D\cap
R_{\varepsilon}}\frac{\rho^n(x)}{Q(x)}\,dm(x)
\end{equation}
для каждого кольца $A(\varepsilon, \varepsilon_0, x_0),$
$\varepsilon_0\in(0,d_0),$ $d_0=\sup\limits_{z\in D}|z-z_0|,$
где $\Sigma_{\varepsilon}$ обозначает семейство всех пересечений
сфер $S(x_0, r)$ с областью $D,$ $r\in (0, \varepsilon_0).$ Примеры
таких отображений несложно указать (см. теорему \ref{thOS4.1} ниже).

Отметим, что выражения $"$почти всех кривых$"$ и $"$почти всех
по\-вер\-х\-но\-с\-тей$"$ в отдельных случаях могут иметь две
различные интерпретации (в частности, если речь идёт о семействе
сфер, то $"$почти всех$"$ может пониматься как относительно
множества значений $r,$ так и конформного модуля семейства сфер,
рассматриваемого как частный случай семейства поверхностей).
Следующее утверждение вносит некоторую ясность между указанными
интерпретациями и может быть установлено полностью по аналогии с
\cite[лемма~9.1]{MRSY}.

\medskip
\begin{lemma}\label{lem8.2.11}{\sl\, Пусть $x_0\in D.$ Если некоторое
свойство $P$ имеет место для почти всех сфер $D(x_0, r):=S(x_0,
r)\cap D,$ где $"$почти всех$"$ понимается в смысле модуля семейств
поверхностей, то $P$ также имеет место для почти всех сфер $D(x_0,
r)$ относительно линейной меры Лебега по параметру $r\in {\Bbb R }.$
Обратно, пусть $P$ имеет место для почти всех сфер $D(x_0,
r):=S(x_0, r)\cap D$ относительно линейной меры Лебега по $r\in
{\Bbb R},$ тогда $P$ также имеет место для почти всех поверхностей
$D(x_0, r):=S(x_0, r)\cap D$ в смысле модуля семейств
поверхностей.}\end{lemma}

\medskip
Следующее утверждение облегчает проверку бесконечной серии
неравенств в (\ref{eq1A}) и может быть установлено аналогично
доказательству \cite[теорема~9.2]{MRSY}.

\medskip
\begin{lemma}\label{lem4}{\sl\,
Пусть $D,$  $D^{\,\prime}\subset\overline{{\Bbb R}^n},$
$x_0\in\overline{D}\setminus\{\infty\}$ и $Q:D\rightarrow(0,\infty)$
-- измеримая по Лебегу функция. Отображение $f:D\rightarrow
D^{\,\prime}$ является нижним $Q$-отображением в точке $x_0$ тогда и
только тогда, когда
%
$$M(f(\Sigma_{\varepsilon}))\ge\int\limits_{\varepsilon}^{\varepsilon_0}
\frac{dr}{||\,Q||_{n-1}(r)}\quad\forall\
\varepsilon\in(0,\varepsilon_0)\,,\ \varepsilon_0\in(0,d_0)\,,$$
%
где, как и выше, $\Sigma_{\varepsilon}$ обозначает семейство всех
пересечений сфер $S(x_0, r)$ с областью $D,$ $r\in (0,
\varepsilon_0),$
$$
\Vert
Q\Vert_{n-1}(r)=\left(\int\limits_{D(x_0,r)}Q^{n-1}(x)\,d{\mathcal{A}}\right)^{\frac{1}{n-1}}$$
-- $L_{n-1}$-норма функции $Q$ над сферой $D(x_0,r)=\{x\in D:
|x-x_0|=r\}=D\cap S(x_0,r)$.}
\end{lemma}

Напомним, что {\it конденсатором} называют пару
$E=\left(A,\,C\right),$ где $A$ -- открытое множество в ${\Bbb
R}^n,$ а $C$ -- компактное подмножество $A.$ {\it Ёмкостью}
конденсатора $E$ называется следующая величина:
\begin{equation}\label{eq1.1AB} {\rm cap}\,E={\rm
cap}\,\left(A,\,C\right)= \inf\limits_{u\in W_0(E)}\,\,\int\limits_A
|\nabla u(x)|^n\,\,dm(x)\,,
\end{equation}
где $W_0(E)=W_0\left(A,\,C\right)$ -- семейство неотрицательных
непрерывных функций $u:A\rightarrow{\Bbb R}$ с компактным носителем
в $A,$ таких что $u(x)\ge 1$ при $x\in C$ и $u\in ACL.$
В формуле (\ref{eq1.1AB}), как обычно, $|\nabla
u|={\left(\sum\limits_{i=1}^n\,{\left(\partial_i u\right)}^2
\right)}^{1/2}.$

Следующее утверждение имеет важное значение для доказательства
многих результатов настоящей работы (см. \cite[предложение~10.2,
гл.~II]{Ri}).

\medskip
\begin{proposition}\label{pr1*!}{\,\sl Пусть $E=(A,\,C)$ --
произвольный конденсатор в ${\Bbb R}^n$ и пусть $\Gamma_E$ --
семейство всех кривых вида $\gamma:[a,\,b)\rightarrow A$ таких, что
$\gamma(a)\in C$ и $|\gamma|\cap\left(A\setminus
F\right)\ne\varnothing$ для произвольного компакта $F\subset A.$
Тогда
%
$${\rm cap}\,E=M(\Gamma_E)\,.$$
%
}
\end{proposition}
В дальнейшем всюду символом $\Gamma(E,F,D)$ мы обозначаем семейство
всех кривых $\gamma:[a,b]\rightarrow\overline{{\Bbb R}^n},$ которые
соединяют $E$ и $F$ в $D,$ т.е. $\gamma(a)\in E,\,\gamma(b)\in F$ и
$\gamma(t)\in D$ при $t\in(a,\,b).$ Для доказательства основных
результатов работы также существенно используются так называемые
кольцевые $Q$-отображения, определение которых приведено ниже (см.,
напр., \cite[гл.~7]{MRSY}, см. также \cite{GRSY} и \cite{BGMV}).
Говорят, что $f:D\rightarrow \overline{{\Bbb R}^n}$ является {\it
кольцевым $Q$-отоб\-ра\-же\-нием в точке $x_0\,\in\,D,$} если
соотношение
%
$$M\left(f\left(\Gamma\left(S_1,\,S_2,\,A\right)\right)\right)\ \le
\int\limits_{A} Q(x)\cdot \eta^n(|x-x_0|)\ dm(x)$$ 
выполнено для любого кольца $A=A(r_1,r_2, x_0),$\, $0<r_1<r_2<
r_0:={\rm dist\,}(x_0, \partial D),$ и для каждой измеримой функции
$\eta : (r_1,r_2)\rightarrow [0,\infty ]\,$ такой, что
%
%
%
$\int\limits_{r_1}^{r_2}\eta(r)dr\ge 1.$
Отметим, что кольцевые $Q$-отображения являются обобщением
квазиконформных отображений и отображений с ограниченным искажением
(см., напр., \cite{GRSY}, \cite{Va}, \cite{Ri}, \cite{BGMV},
\cite{Pol}, \cite{Re$_1*$}, \cite{Re} и \cite{Vu}). В частности,
ввиду неравенства Е.~Полецкого отображения с ограниченным искажением
являются кольцевыми $Q$-отображениями с некоторой ограниченной
функцией $Q$ (см. \cite[теорема~1]{Pol}). Следующее утверждение было
доказано и опубликовано автором данной работы несколько ранее и
может быть найдено в \cite[теорема~1]{Sev}.

\medskip
\begin{lemma}\label{lem2} {\sl\, Пусть $Q:D\rightarrow [0, \infty]$ -- измеримая по Лебегу
функция, $Q\in L_{loc}^1(D).$ Открытое дискретное отображение
$f:D\rightarrow \overline{{\Bbb R}^n}$ является кольцевым
$Q$-отображением в точке $x_0\in D$ тогда и только тогда, когда для
произвольных $0<r_1<r_2< {\rm dist} \, (x_0,\partial D)$ и
произвольного конденсатора $E=\left(B(x_0, r_2), \overline{B(x_0,
r_1)}\right)$ ёмкость конденсатора $$f(E):=\left(f(B(x_0, r_2)),
f\left(\overline{B(x_0, r_1)}\right)\right)$$ удовлетворяет условию
%
$${\rm cap}\, f(E)\le\frac{\omega_{n-1}}{I^{n-1}}\,,$$
%
где $I=I(x_0,r_1,r_2)$ задаётся соотношением
%
$$I=I(x_0,r_1,r_2)=\int\limits_{r_1}^{r_2}\
\frac{dr}{rq_{x_0}^{\frac{1}{n-1}}(r)}\,.$$
}
\end{lemma}

Следующие важные сведения, касающиеся ёмкости пары множеств
относительно области, могут быть найдены в работе В.~Цимера
\cite{Zi}. Пусть $G$ -- ограниченная область в ${\Bbb R}^n$ и $C_{0}
, C_{1}$ -- непересекающиеся компактные множества, лежащие в
замыкании $G.$ Полагаем  $R=G \setminus (C_{0} \cup C_{1})$ и
$R^{\,*}=R \cup C_{0}\cup C_{1}.$ {\it Конформной ёмкостью пары
$C_{0}, C_{1}$ относительно замыкания $G$} называется величина
$$C[G, C_{0}, C_{1}] = \inf \int\limits_{R} \vert \nabla u
\vert^{n}\ dm(x),$$
где точная нижняя грань берётся по всем функциям $u,$ непрерывным в
$R^{\,*},$ $u\in ACL(R),$ таким что $u=1$ на $C_{1}$ и $u=0$ на
$C_{0}.$ Указанные функции будем называть {\it допустимыми} для
величины $C [G, C_{0}, C_{1}].$ Мы будем говорить, что  {\it
множество $\sigma \subset {\Bbb R}^n$ разделяет $C_{0}$ и $C_{1}$ в
$R^{\,*}$} если $\sigma \cap R$ замкнуто в $R$ и найдутся
непересекающиеся множества $A$ и $B,$ являющиеся открытыми в
$R^{\,*} \setminus \sigma,$ такие что $R^{\,*} \setminus \sigma =
A\cup B,$ $C_{0}\subset A$ и $C_{1} \subset B.$ Пусть $\Sigma$
обозначает класс всех множеств, разделяющих $C_{0}$ и $C_{1}$ в
$R^{\,*}.$ Для числа $n^{\prime} = n/(n-1)$ определим величину
\begin{equation}\label{eq13.4.12}
\widetilde{M_{n^{\prime}}}(\Sigma)=\inf\limits_{\rho\in
\widetilde{\rm adm} \Sigma} \int\limits_{{\Bbb
R}^n}\rho^{\,n^{\prime}}dm(x)
\end{equation}
где запись $\rho\in \widetilde{\rm adm}\,\Sigma$ означает, что
$\rho$ -- неотрицательная борелевская функция в ${\Bbb R}^n$ такая,
что
\begin{equation} \label{eq13.4.13}
\int\limits_{\sigma \cap R}\rho d{\mathcal H}^{n-1} \ge
1\quad\forall\, \sigma \in \Sigma\,. \end{equation}
Заметим, что согласно результата Цимера,
\begin{equation}\label{eq3}
\widetilde{M_{n^{\,\prime}}}(\Sigma)=C [G , C_{0} ,
C_{1}]^{\,-1/(n-1)}\,,
\end{equation}
см. \cite[теорема~3.13]{Zi}. Заметим также, что согласно результата
Хессе
\begin{equation}\label{eq4}
M(\Gamma(E, F, D))= C[D, E, F]\,,
\end{equation}
как только $(E \cap F)\cap
\partial D = \varnothing,$
см. \cite[теорема~5.5]{Hes}.

\medskip
Сформулируем и докажем следующие утверждения.

\medskip
\begin{lemma}\label{lem1} {\sl Пусть $x_0\in D$ и $Q:D\rightarrow
[0,\infty]$ -- локально интегрируемая в степени $n-1$ в $D$ функция.
Если $f:D\rightarrow \overline{{\Bbb R}^n}$ -- локальный нижний
$Q$-гомеоморфизм в точке $x_0,$ то $f$ является кольцевым
$Q^{\,*}$-отображением в этой же точке при
$Q^{\,*}=Q^{n-1}.$}\end{lemma}

\begin{proof} Пусть $x_0\in D,$ $0<r_1<r_2<{\rm dist\,}(x_0, \partial D).$
Без ограничения общности, мы можем считать, что $f(x_0)\ne \infty.$
Согласно лемме \ref{lem2} достаточно установить, что
$${\rm cap}\,
f(E)\le \frac{\omega_{n-1}}{I^{*\,n-1}}$$
где $E$ -- конденсатор вида $E=(B(x_0, r_2), \overline{B(x_0,
r_1)}),$ $\omega_{n-1}$ -- площадь единичной сферы в ${\Bbb R}^n,$
$q^{\,*}_{x_0}(r)$ -- среднее значение функции $Q^{n-1}(x)$ над
сферой $|x-x_0|=r$ и $I^{\,*}=I^{\,*}(x_0,
r_1,r_2)=\int\limits_{r_1}^{r_2}\
\frac{dr}{rq^{\,*\,\frac{1}{n-1}}_{x_0}(r)}.$ Зафиксируем
$\varepsilon\in (r_1, r_2)$ и рассмотрим шар $B(x_0, \varepsilon).$
Полагаем $C_0=\partial f(B(x_0, r_2)),$ $C_1=f(\overline{B(x_0,
r_1)}),$ $\sigma=\partial f(B(x_0, \varepsilon)).$ Поскольку
$\overline{B(x_0, r_2)}$ -- компакт в $D,$ найдётся шар $B(x_0, R)$
такой, что $\overline{f(B(x_0, r_2))}\subset B(x_0, R).$ Полагаем
$G:=B(x_0, R).$

\medskip
Поскольку $f$ -- локальный гомеоморфизм, $\overline{f(B(x_0, r_1))}$
-- компактное подмножество $f(B(x_0, \varepsilon))$ также, как
$\overline{f(B(x_0, \varepsilon))}$ -- компактное подмножество
$f(B(x_0, r_2)).$ В частности, $\overline{f(B(x_0, r_1))}\cap
\partial f(B(x_0, \varepsilon))=\varnothing.$ Пусть, как и выше, $R=G
\setminus (C_{0} \cup C_{1})$ и $R^{\,*} = R \cup C_{0}\cup C_{1},$
тогда $R^{\,*}:=G.$ Заметим, что $\sigma$ разделяет $C_0$ и $C_1$ в
$R^{\,*}=G.$ Действительно, множество $\sigma \cap R$ замкнуто в
$R,$ кроме того, пусть $A:=G\setminus \overline{f(B(x_0,
\varepsilon))}$ и $B= f(B(x_0, \varepsilon)),$ тогда $A$ и $B$
открыты в $G\setminus \sigma,$ $C_0\subset A,$ $C_1\subset B$ и
$G\setminus \sigma=A\cup B.$

\medskip
Пусть $\Sigma$ -- семейство всех множеств, отделяющих $C_0$ от $C_1$
в $G.$ Поскольку для открытых отображений $\partial f(O)\subset
f(\partial O),$ где $O$ -- компактная подобласть $D,$ мы получим:
$\partial f(B(x_0, r))\subset f(\partial B(x_0, r)),$ $r\in (0, {\rm
dist\,}(x_0, \partial D)).$

\medskip
Пусть $\rho^{n-1}\in \widetilde{{\rm
adm}}\bigcup\limits_{r=r_1}^{r_2}
\partial f(B(x_0, r))$ в смысле соотношения (\ref{eq13.4.13}), тогда также
$\rho\in {\rm adm}\bigcup\limits_{r=r_1}^{r_2}
\partial f(B(x_0, r))$ в смысле соотношения (\ref{eq8.2.6}).
Поскольку (ввиду открытости отображения $f$) имеет место включение
$\partial f(B(x_0, r))\subset f(S(x_0, r)),$ мы получим, что
$\rho\in {\rm adm}\bigcup\limits_{r=r_1}^{r_2} f(S(x_0, r))$ и,
следовательно, ввиду (\ref{eq13.4.12}) будем иметь
$$\widetilde{M_{n^{\prime}}}(\Sigma)\ge
\widetilde{M_{n^{\prime}}}\left(\bigcup\limits_{r=r_1}^{r_2}
\partial f(B(x_0, r))\right)\ge M\left(\bigcup\limits_{r=r_1}^{r_2}
\partial f(B(x_0, r))\right)\ge
$$
\begin{equation}\label{eq5}
\ge M\left(\bigcup\limits_{r=r_1}^{r_2} f(S(x_0, r))\right)\,.
\end{equation}
Однако, ввиду (\ref{eq3}) и (\ref{eq4}),
\begin{equation}\label{eq6}
\widetilde{M_{n^{\prime}}}(\Sigma)=\frac{1}{(M(\Gamma(C_0, C_1,
G)))^{1/(n-1)}}\,.
\end{equation}
Пусть $\Gamma_{f(E)}$ -- семейство всех кривых для конденсатора
$f(E)$ в обозначениях предложения \ref{pr1*!}. Пусть также
$\Gamma^{\,*}_{f(E)}$ обозначает семейство всех спрямляемых кривых
семейства $\Gamma_{f(E)},$ тогда заметим, что семейства
$\Gamma^{*}_{f(E)}$ и $\Gamma (C_0, C_1, G)$ имеют одинаковые
семейства допустимых метрик $\rho$ и, значит,
$$M(\Gamma_{f(E)})=M(\Gamma(C_0, C_1, G))\,.$$ Из (\ref{eq6}) и
предложения \ref{pr1*!} мы получим, что
\begin{equation}\label{eq7}
\widetilde{M}^{n-1}(\Sigma)=\frac{1}{{\rm cap\,}f(E)}\,.
\end{equation}
Окончательно, из (\ref{eq5}) и (\ref{eq7}) мы получаем неравенство
\begin{equation}\label{eq8}
{\rm cap\,}f(E) \le \frac{1}{M\left(\bigcup\limits_{r=r_1}^{r_2}
f(S(x_0, r))\right)^{n-1}}\,.
\end{equation}
По лемме \ref{lem4} и из (\ref{eq8}) мы получим, что
$$
{\rm cap\,}f(E) \le \frac{1}{\left(\int\limits_{r_1}^{r_2}
\frac{dr}{\Vert
\,Q\Vert_{n-1}(r)}\right)^{n-1}}=\frac{\omega_{n-1}}{I^{*\,n-1}}\,,
$$
что и доказывает утверждение леммы \ref{lem1}.
\end{proof}$\Box$

\medskip
Напомним, что отображение $f:X\rightarrow Y$ между пространствами с
мерами $(X, \Sigma, \mu)$ и $(Y, \Sigma^{\,\prime}, \mu^{\,\prime})$
обладает {\it $N$-свой\-с\-т\-вом} (Лузина), если из условия
$\mu(S)=0$ следует, что $\mu^{\,\prime}(f(S))=0.$ Следующее
вспомогательное утверждение получено в работе \cite{KRSS} (см.
теорема 1 и следствие 2).

\medskip
\begin{proposition}\label{pr1}
{\sl\, Пусть $D$ -- область в ${\Bbb R}^n,$ $n\ge 3,$
$\varphi:(0,\infty)\rightarrow (0,\infty)$ -- неубывающая функция,
удовлетворяющая условию (\ref{eqOS3.0a}). Тогда:

1) Если $f:D\rightarrow{\Bbb R}^n$ -- непрерывное открытое
отображение класса $W^{1,\varphi}_{loc}(D),$ то $f$ имеет почти
всюду полный дифференциал в $D;$

2) Любое непрерывное отображение $f\in W^{1,\varphi}_{loc}$ обладает
$N$-свойством относительно $(n-1)$-мерной меры Хаусдорфа, более
того, локально абсолютно непрерывно на почти всех сферах $S(x_0, r)$
с центром в заданной предписанной точке $x_0\in{\Bbb R}^n$. Кроме
того, на почти всех таких сферах $S(x_0, r)$ выполнено условие
${\mathcal H}^{n-1}(f(E))=0,$ как только $|\nabla f|=0$ на множестве
$E\subset S(x_0, r).$ (Здесь $"$почти всех$"$ понимается
относительно линейной меры Лебега по параметру $r$).}

\end{proposition}

Имеет место следующее утверждение.

\medskip
\begin{theorem}{}\label{thOS4.1} {\sl Пусть $D$ -- область в ${\Bbb R}^n,$
$n\ge 3,$ $\varphi:(0,\infty)\rightarrow (0,\infty)$ -- неубывающая
функция, удовлетворяющая условию (\ref{eqOS3.0a}).
Если $n\ge 3,$ то каждый локальный гомеоморфизм $f:D\rightarrow
{\Bbb R}^n$ с конечным искажением класса $W^{1,\varphi}_{loc}$
такой, что $N(f, D)<\infty,$ является нижним $Q$-отображением в
каждой точке $x_0\in\overline{D}$ при $Q(x)=N(f, D)\cdot K_O(x, f),$
где внешняя дилатация $K_O(x, f)$ отображения $f$ в точке $x$
определена соотношением (\ref{eq0.1.1A}), а кратность $N(f, D)$
определена соотношением (\ref{eq1.7A}).}
\end{theorem}

\medskip
\begin{proof}
Заметим, что $f$ дифференцируемо почти всюду ввиду предложения
\ref{pr1}. Пусть $B$ -- борелево множество всех точек $x\in D,$ в
которых $f$ имеет полный дифференциал $f^{\,\prime}(x)$ и $J(x,
f)\ne 0.$ Применяя теорему Кирсбрауна и свойство единственности
аппроксимативного дифференциала (см. \cite[пункты~2.10.43 и
3.1.2]{Fe}), мы видим, что множество $B$ представляет собой не более
чем счётное объединение борелевских множеств $B_l,$
$l=1,2,\ldots\,,$ таких, что сужения $f_l=f|_{B_l}$ являются
билипшецевыми гомеоморфизмами (см., напр., \cite[пункты~3.2.2, 3.1.4
и 3.1.8]{Fe}). Без ограничения общности, мы можем полагать, что
множества $B_l$ попарно не пересекаются. Обозначим также символом
$B_*$ множество всех точек $x\in D,$ в которых $f$ имеет полный
дифференциал, однако, $f^{\,\prime}(x)=0.$

\medskip
Ввиду построения, множество $B_0:=D\setminus \left(B\bigcup
B_*\right)$ имеет лебегову меру нуль. Следовательно, по
\cite[теорема~9.1]{MRSY}, ${\mathcal H}^{n-1}(B_0\cap S_r)=0$ для
почти всех сфер $S_r:=S(x_0,r)$ с центром в точке
$x_0\in\overline{D},$ где $"$почти всех$"$ следует понимать в смысле
конформного модуля семейств поверхностей. По лемме \ref{lem8.2.11}
также ${\mathcal H}^{n-1}(B_0\cap S_r)=0$ при почти всех $r\in {\Bbb
R}.$

\medskip
Ввиду предложения \ref{pr1} из условия ${\mathcal H}^{n-1}(B_0\cap
S_r)=0$ для почти всех $r\in {\Bbb R}$ вытекает, что ${\mathcal
H}^{n-1}(f(B_0\cap S_r))=0$ для почти всех $r\in {\Bbb R}.$ По этому
предложению также ${\mathcal H}^{n-1}(f(B_*\cap S_r))=0,$ поскольку
$f$--отображение с конечным искажением и, значит, $\nabla f=0$ почти
всюду, где $J(x, f)=0.$

\medskip
Пусть $\Gamma$ -- семейство всех пересечений сфер $S_r,$
$r\in(\varepsilon,\varepsilon_0),$
$\varepsilon_0<d_0=\sup\limits_{x\in D}\,|x-x_0|,$  с областью $D.$
Для заданной функции $\rho_*\in{\rm adm}\,f(\Gamma),$
$\rho_*\equiv0$ вне $f(D),$ полагаем $\rho\equiv 0$ вне $D$ и на
$B_0,$
$$\rho(x)\ \colon=\ \rho_*(f(x))\Vert f^{\,\prime}(x)\Vert
\qquad\text{при}\ x\in D\setminus B_0\,.$$
Пусть $D_{r}^{\,*}\in f(\Gamma),$ $D_{r}^{\,*}=f(D\cap S_r).$
Заметим, что
$$D_{r}^{\,*}=\bigcup\limits_{i=0}^{\infty} f(S_r\cap
B_i)\bigcup f(S_r\cap B_*)\,,$$
и, следовательно, для почти всех $r\in (0, \varepsilon_0)$
\begin{equation}\label{eq10B}
1\le \int\limits_{D^{\,*}_r}\rho^{n-1}_*(y)d{\mathcal A_*} \le
\sum\limits_{i=0}^{\infty} \int \limits_{f(S_r\cap B_i)}
\rho^{n-1}_*(y)N (y, S_r\cap B_i)d{\mathcal H}^{n-1}y +
\end{equation}
$$+\int\limits_{f(S_r\cap B_*)} \rho^{n-1}_*(y)N (y, S_r\cap B_*
)d{\mathcal H}^{n-1}y\,.$$ Учитывая доказанное выше, из
(\ref{eq10B}) мы получаем, что
\begin{equation}\label{eq11B}
1\le \int\limits_{D^{\,*}_r}\rho^{n-1}_*(y)d{\mathcal A_*} \le
\sum\limits_{i=1}^{\infty} \int \limits_{f(S_r\cap B_i)}
\rho^{n-1}_*(y)N (y, S_r\cap B_i)d{\mathcal H}^{n-1}y
\end{equation}
для почти всех $r\in (0, \varepsilon_0).$
Рассуждая покусочно на $B_i,$ $i=1,2,\ldots,$ ввиду \cite[1.7.6 и
теорема~3.2.5]{Fe} мы получаем, что
$$\int\limits_{B_i\cap S_r}\rho^{n-1}\,d{\mathcal A}=
\int\limits_{B_i\cap S_r}\rho_*^{n-1}(f(x))\Vert
f^{\,\prime}(x)\Vert^{n-1}\,d{\mathcal A}=$$
$$=\int\limits_{B_i\cap S_r}\rho_*^{n-1}(f(x))\cdot\frac{\Vert
f^{\,\prime}(x)\Vert^{n-1}}{\frac{d{\mathcal A_*}}{d{\mathcal
A}}}\cdot \frac{d{\mathcal A_*}}{d{\mathcal A}}\,d{\mathcal A}\ge
\int\limits_{B_i\cap S_r}\rho_*^{n-1}(f(x))\cdot \frac{d{\mathcal
A_*}}{d{\mathcal A}}\,d{\mathcal A}=$$
\begin{equation}\label{eq12B}
=\int\limits_{f(B_i\cap S_r)}\rho_{*}^{n-1}\,N(y, S_r\cap
B_i)d{\mathcal H}^{n-1}y
\end{equation} для почти всех $r\in (0, \varepsilon_0).$
Из (\ref{eq11B}) и (\ref{eq12B}) вытекает, что
$\rho\in{\rm{ext\,adm}}\,\Gamma.$

Замена переменных на каждом $B_l,$ $l=1,2,\ldots\,,$ см., напр.,
\cite[теорема~3.2.5]{Fe}, и свойство счётной аддитивности интеграла
приводят к оценке
$$\int\limits_{D}\frac{\rho^n(x)}{K_O(x, f)}\,dm(x)\le
\int\limits_{f(D)}N(f, D)\cdot \rho^{\,n}_*(y)\, dm(y)\,,$$ что и
завершает доказательство.
\end{proof}$\Box$

\medskip
{\bf 3. Доказательство основных результатов}. Следующий результат
даже в несколько более общем виде доказан в работе
\cite[теорема~1]{GolSev}.

\medskip
\begin{proposition}\label{pr2}
{\sl Пусть $f:{\Bbb B}^n\rightarrow {\Bbb R}^n,$ $n\ge 3,$ --
локальный кольцевой $Q$-го\-ме\-о\-мор\-физм в точке $x_0=0,$ такой
что $Q\in L_{loc}^{1}({\Bbb B}^n).$ Если выполнено хотя бы одно из
условий: (\ref{eq9}), либо $Q\in FMO(0),$
то отображение $f$ инъективно в некотором шаре $B\left(0, \delta(n,
Q)\right),$ где $\delta$ -- положительное число, зависящее только от
$n$ и функции $Q.$}
\end{proposition}

\medskip
{\it Доказательство теорем \ref{th1} и \ref{th2}.} Ввиду леммы
\ref{lem1} и теоремы \ref{thOS4.1} каждое отображение
$f\in\frak{R}_{\varphi, Q, N, \delta}$ является кольцевым локальным
$Q\cdot N^{n-1}$-гомеоморфизмом в каждой точке области $D.$ Ввиду
предложения \ref{pr2} каждая точка имеет окрестность, радиус которой
зависит только от $n,$ $Q$ и числа $N,$ в которой $f$ является
гомеоморфизмом. Необходимое заключение вытекает, в таком случае, из
результатов, полученных для $Q$-гомеоморфизмов ранее (см.
\cite[теоремы~4.1--4.2]{RS}). $\Box$

\medskip
{\it Доказательство теоремы \ref{th3}.} Ввиду леммы \ref{lem1} и
теоремы \ref{thOS4.1} каждое отображение $f_m$ из формулировки
теоремы \ref{thOS4.1} является кольцевым локальным $Q\cdot
N^{n-1}$-гомеоморфизмом в каждой точке области $D.$ Пусть $G$ --
произвольная компактная подобласть области $D.$ Покроем
$\overline{G}$ (возможно) бесконечным числом шаров с центрами во
всех точках $\overline{G}$ и имеющих радиусы, соответствующие
утверждению предложения \ref{pr2} (т.е., радиусы шаров, зависящие
только только от $n,$ $Q$ и числа $N,$ в которых отображение $f_m$
является инъективным). По лемме Гейне--Бореля--Лебега можно
ограничиться лишь конечным набором таких шаров. Заметим, что ввиду
\cite[теоремы~4.1--4.2]{RSS} предельное отображение $f$ является
гомеоморфизмом либо постоянной в каждом таком шаре. Пусть $f$
постоянна в каком либо шаре $B_1,$ тогда из связности $\overline{G}$
вытекает, что этот шар пересекается по крайней мере ещё с одним
шаром $B_2,$ т.о., $f\equiv const$ в $B_2$ и т.д. Продолжая этот
процесс, мы приходим к заключению, что в этом случае $f\equiv const$
в $\overline{G}.$ Устроив исчерпание области $D$ компактными
подобластями $G,$ мы можем заключить, что $f$ постоянно в $D.$ Пусть
теперь $f$ -- локальный гомеоморфизм в каждом шаре указанного выше
покрытия, тогда, очевидно, $f$ -- локальный гомеоморфизм в $G$ и,
значит, $f$ -- локальный гомеоморфизм в $D,$ что и требовалось
установить. $\Box$

\medskip
КОНТАКТНАЯ ИНФОРМАЦИЯ

\medskip
\noindent{{\bf Евгений Александрович Севостьянов} \\ Институт
прикладной математики и механики НАН Украины \\
83 114 Украина, г. Донецк, ул. Розы Люксембург, д. 74, \\
тел. +38 (066) 959 50 34 (моб.), +38 (062) 311 01 45 (раб.), e-mail:
brusin2006@rambler.ru, esevostyanov2009@mail.ru}

\end{document}